\documentclass[11pt,a4paper,reqno]{amsart}
\usepackage[T1]{fontenc}
\usepackage[utf8]{inputenc}
\usepackage{lmodern}
\usepackage{amsmath,amssymb,amsthm,mathtools}
\usepackage{enumitem}
\usepackage{needspace}
\usepackage{microtype}
\usepackage{bm}
\usepackage[colorlinks=true,linkcolor=blue,citecolor=blue,urlcolor=blue]{hyperref}
\numberwithin{equation}{section}

\newtheorem{thm}{Theorem}[section]
\newtheorem{prop}[thm]{Proposition}
\newtheorem{lem}[thm]{Lemma}
\newtheorem{cor}[thm]{Corollary}
\theoremstyle{definition}

\newtheorem{question}[thm]{Question}
\newtheorem{rem}[thm]{Remark}
\newtheorem{exa}[thm]{Example}

\newcommand{\R}{\mathbb R}

\newcommand{\Q}{\mathbb Q}

\title[]{
On regularity-free characterization of bisymmetric means}

\author{Gergely Kiss}
\address{Gergely Kiss \newline
HUN-REN Alfr\'ed R\'enyi Institute of Mathematics, Analysis Department, Reáltanoda u. 13-15, 1053, Budapest, Hungary, \newline
and \newline
Corvinus University of Budapest, F\"ov\'am t\'er 13-15, 1093, Budapest, Hungary}
\email{kigergo57@gmail.com}
\author{Kazuki Okamura}
\address{Kazuki Okamura \newline Department of Mathematics, Faculty of Science, Shizuoka University, 836 Ohya, Suruga-ku, Shizuoka 422-8529, Japan}
\email{okamura.kazuki@shizuoka.ac.jp}

\date{\today}

\begin{document}

\begin{abstract}
We prove that every reflexive, bisymmetric binary operation on a non-degenerate real interval that is strictly increasing in each variable is automatically continuous and is a weighted quasi-arithmetic mean. This settles the reflexive, strictly increasing case of the automatic-continuity problem without assuming symmetry.
The proof constructs a symmetric mean as the unique fixed point of a natural increasing operation and then applies the known symmetric automatic-continuity theorem.
Maksa's method, combined with an interpolation argument on a countable dense set, yields an exact characterization of operations of arbitrary finite arity, with uniquely determined positive weights.
\end{abstract}

\maketitle

\section{Introduction}

A weighted quasi-arithmetic mean on a non-degenerate real interval $I$ has the form
\[ M(x,y)=\varphi^{-1}\left(p\varphi(x)+(1-p)\varphi(y)\right),
 \qquad 0<p<1,\]
where $\varphi : I\to\R$ is continuous and strictly increasing.
The characterization of these means by reflexivity, 
bisymmetry, strict increase and continuity is classical; see Acz\'el~\cite{Aczel1948} and Acz\'el--Dhombres~\cite{AczelDhombres1989}.
The question addressed here is whether continuity follows from the other three assumptions.

The aim of this paper is twofold.
First, we prove that symmetry can be omitted from the regularity-free characterization of reflexive, bisymmetric, strictly increasing operations.
We also review the history of continuity-free characterization of bisymmetric operations.
The binary result below answers the reflexive and strictly increasing case of the continuity questions discussed by Burai, Kiss and Szokol~\cite{BuraiKissSzokol2021,BuraiKissSzokol2023}. 
In parallel, Burai, Kiss and Szokol  \cite{BuraiKissSzokol2026} are developing a constructive proof for algebraic values of the weight parameter $p$ in a separate work in preparation.

Let $I$ be a non-degenerate real interval and let $M:I^2\to I$ be a binary operation.
We call $M$
\begin{itemize}
\item \emph{reflexive} if $M(x,x)=x$ for every $x\in I$,
\item \emph{symmetric} if $M(x,y)=M(y,x)$ for every $x,y\in I$,
\item \emph{bisymmetric} if
\[ M(M(x,y),M(u,v))=M(M(x,u),M(y,v)),  \qquad x,y,u,v\in I.\]
\end{itemize}
We say that $M$ is \emph{strictly increasing in each variable}
if, for every $y\in I$, both maps $x\mapsto M(x,y)$ and $x\mapsto M(y,x)$ are strictly increasing.

We call $M$ a \emph{mean} if
\[ \min\{x,y\}\leq M(x,y)\leq\max\{x,y\},  \qquad x,y\in I,\]
and a \emph{strict mean} if both inequalities are strict whenever $x \ne y$.
Let $F: I^2\to I$ be an operation which is strictly increasing in each variable.
Then
\[    F \textrm{ is a mean }\quad \Longleftrightarrow \quad F \textrm{ is a strict mean}\quad \Longleftrightarrow \quad F \textrm{ is reflexive}.\]

Indeed, every mean is reflexive, while reflexivity and strict increase give that for $x<y$,
\[
 \begin{aligned}
 x&=M(x,x)<M(x,y)<M(y,y)=y,\\
 x&=M(x,x)<M(y,x)<M(y,y)=y.
 \end{aligned}
\]
However, there exist strict means which are nondecreasing in each variable but discontinuous (see Example~\ref{exa:strict-mean}). The operation in that example is reflexive, symmetric, nondecreasing and bisymmetric, yet it is not a weighted quasi-arithmetic mean.
On the other hand, a continuous bisymmetric operation is a strict mean if and only if it is reflexive and strictly increasing in each variable (see Proposition~\ref{prop:continuous-strict-mean} and \cite{FodorMarichal1997}).
This distinction motivates the search for a regularity-free characterization of 
bisymmetric strict means that are strictly increasing in each variable.
Our first main result in this direction is the following.

\Needspace{10\baselineskip}
\begin{thm}\label{thm:main}
Let $I$ be a non-degenerate real interval and let $M:I^2\to I$ be reflexive, bisymmetric and strictly increasing in each variable.
Then there exist a unique number $p\in(0,1)$ and a continuous strictly increasing function $\varphi:I\to\R$, unique up to a positive affine change, such that
\begin{equation}\label{eq:main-representation}
 M(x,y)=\varphi^{-1}
 \left(p\varphi(x)+(1-p)\varphi(y)\right),
 \qquad x,y\in I.
\end{equation}
In particular, $M$ is continuous.
\end{thm}

For $x,y,t\in I$, let
\[\Phi_{x,y}(t)\coloneqq M(M(t,x),M(y,t)).\]
Acz\'el \cite[pp.~396--398]{Aczel1948} used this map to identify the limiting symmetric mean under continuity.
Here we show that $\Phi_{x,y}$ has a unique fixed point $m(x,y)$ directly from order and bisymmetry.

The resulting operation $m$ is reflexive, symmetric, bisymmetric and strictly increasing in each variable, and satisfies
\begin{equation}\label{eq:intro-intertwining}
 M(m(x,y),m(z,w))=m(M(x,z),M(y,w)).
\end{equation}
The symmetric automatic-continuity theorem of Burai, Kiss and Szokol~\cite{BuraiKissSzokol2021} transforms $m$ into the arithmetic midpoint.
In the same coordinates, \eqref{eq:intro-intertwining} becomes Jensen's midpoint equation for the transformed operation $M$.
Monotonicity and reflexivity then give the weighted representation.

The binary theorem also gives the complete characterization in every finite arity.
Here bisymmetry means that evaluation by rows and by columns gives the same result for every square array.
Its precise form is stated in~\eqref{eq:nary-bisymmetry}.

\begin{thm}\label{thm:nary-main}
Let $n\geq2$ and let $I$ be a non-degenerate real interval.
An operation $F:I^n\to I$ is reflexive, bisymmetric and strictly increasing in each variable if and only if there exist positive numbers $\lambda_1,\ldots,\lambda_n$ with
\[ \lambda_1+\cdots+\lambda_n=1\]
and a continuous strictly increasing function $\varphi:I\to\R$ such that
\begin{equation}\label{eq:nary-form}
 F(x_1,\ldots,x_n) =\varphi^{-1}\!\left(\sum_{i=1}^n\lambda_i\varphi(x_i)\right),
 \qquad x_1,\ldots,x_n\in I.
\end{equation}
The weights are uniquely determined, and the generator $\varphi$ is unique up to a positive affine change.
In particular, $F$ is continuous.
\end{thm}

Section~\ref{sec:history} places these results in their historical context.
Section~\ref{sec:symmetrization} proves Theorem~\ref{thm:main}, and Section~\ref{sec:nary} proves Theorem~\ref{thm:nary-main} using a binary retract and an elementary interpolation argument.
As an application of our results, Section \ref{sec:strong-decomposability} deals with the strongly decomposable families considered by Marichal \cite{Marichal2000} without regularity.
Section~\ref{sec:questions} discusses weighted families, nonstrict monotonicity and affine representations beyond the reflexive case.

\section{History of the problem}\label{sec:history}

The classical theory starts with the characterizations of quasi-arithmetic means by Kolmogorov~\cite{Kolmogorov1930}, Nagumo~\cite{Nagumo1930} and de~Finetti~\cite{deFinetti1931}.
Their setting is a compatible family $(M_n)_{n\geq1}$ of means, represented by a single continuous strictly increasing generator:
\[M_n(x_1,\ldots,x_n)
=\varphi^{-1}\!\left(\frac1n\sum_{i=1}^n\varphi(x_i)\right).\]
Besides continuity, symmetry, reflexivity and strict increase, the characterizations use a replacement axiom: replacing any block of inputs by repeated copies of its own mean leaves the overall mean unchanged.
Thus compatibility relates operations of different arities.
Recently, Kiss and Shulman~\cite{KissShulman} remove the continuity assumption from this family characterization.

Acz\'el~\cite{Aczel1948} introduced a complementary approach based on a single binary operation.
In the symmetric case, continuity, reflexivity, bisymmetry and strict monotonicity characterize the quasi-arithmetic midpoint.
Dropping symmetry in the continuous, strictly increasing case gives weighted quasi-arithmetic means.
The passage to fixed higher arity and to generalized bisymmetry was developed by Acz\'el--Maksa~\cite{AczelMaksa1996}, Maksa~\cite{Maksa1999,Maksa2002}, and M\"unnich--Maksa--Mokken~\cite{MaksaMokkenMunnich}.
In particular, Maksa's proof in~\cite{Maksa2002} uses a binary retract and a fixed-weight Jensen equation; Section~\ref{sec:nary} follows this reduction after obtaining the missing binary continuity theorem.

Quasi-arithmetic means, bisymmetry and their generalizations have continued to attract substantial attention.
In addition to the preceding work on generalized bisymmetry, representative directions include nonstrict means and axiomatizations without symmetry~\cite{FodorMarichal1997,Marichal2000}, Gauss composition and equations involving weighted means~\cite{DaroczyPales2002,DaroczyMaksa2013}, characterizations of generalized quasi-arithmetic means~\cite{MatkowskiPales2015}, equality and comparison problems~\cite{Pales2011,PalesPasteczka2026}, and invariance~\cite{JarczykJarczyk2018}.

There is also a parallel universal-algebraic strand, less often connected explicitly with the analytic theory of means.
In that language, bisymmetry is usually called \emph{mediality}, reflexivity is \emph{idempotence}, and strict increase implies cancellation.
The affine representations of medial quasigroups in the work of Toyoda~\cite{Toyoda1941} and Bruck~\cite{Bruck1944}, Evans's theory of abstract mean values~\cite{Evans1963}, and Howroyd's work on cancellative medial groupoids~\cite{Howroyd1973} exhibit the algebraic structure behind these identities.
The present fixed-point argument provides a direct order-theoretic route to the representation without requiring the quasigroup division operations.

Continuity is an explicit hypothesis in Acz\'el's original characterization.
Burai, Kiss and Szokol~\cite{BuraiKissSzokol2021} prove that it follows from the remaining assumptions in the symmetric, reflexive binary case.
The following questions for operations on a non-degenerate real interval $I$ can be found in~\cite{BuraiKissSzokol2021,BuraiKissSzokol2023}:
\begin{enumerate}[label=\textup{(\roman*)},leftmargin=*]
\item Is every reflexive, bisymmetric operation $F:I^2\to I$ which is strictly increasing in each variable necessarily continuous?
\item Is every bisymmetric operation $F:I^2\to I$ which  is strictly increasing in each variable necessarily continuous?
\end{enumerate}
The dichotomy theorem of Burai, Kiss and Szokol~\cite{BuraiKissSzokol2023} shows that a reflexive, bisymmetric, strictly increasing binary operation is symmetric and continuous as soon as symmetry holds at one off-diagonal pair.
Theorem~\ref{thm:main} answers~\textup{(i)} affirmatively without that symmetry assumption; Theorem~\ref{thm:nary-main} gives the corresponding result in every finite arity.

Question~\textup{(ii)} has a negative answer in every finite arity.
Kiss~\cite[Theorems~3.4 and~4.1]{Kiss2026} constructs, for each $n\geq2$ and each choice of $\alpha_1,\ldots,\alpha_n>0$ with $\sum_i\alpha_i\ne1$, a bisymmetric operation
\begin{equation*}
 F(x_1,\ldots,x_n) =\psi^{-1}\!\left(\sum_{i=1}^n\alpha_i\psi(x_i)\right)
\end{equation*}
on a suitable real interval, strictly increasing in every variable, for which every one-variable section is discontinuous.
Taking equal coefficients also gives symmetric examples, showing that symmetry alone does not restore automatic continuity without reflexivity.
Here $\psi$ is an increasing bijection onto a set closed under the indicated weighted sums; its image need not be an interval, and continuity of $\psi$ is not assumed.
In the binary case, this reads
\[
F(x,y)=\psi^{-1}\left(\alpha\psi(x)+\beta\psi(y)\right).\]
For positive coefficients, an operation of this form is reflexive exactly when their sum is~$1$.
In that case, Theorem~\ref{thm:nary-main} forces continuity.
Together with the construction, this gives the precise boundary.
A discontinuous example with prescribed positive coefficients exists if and only if their sum differs from~$1$.
This is an existence statement.
Indeed, coefficients whose sum differs from~$1$ also admit continuous examples, for instance with $I=\R$ and $\psi$ the identity.

\section{Canonical fixed-point symmetrization}\label{sec:symmetrization}

Throughout this section, $I$ is a non-degenerate real interval and $M:I^2\to I$ is reflexive, bisymmetric and strictly increasing in each variable.


For $x,y\in I$, define an increasing map
\begin{equation*}
 \Phi_{x,y}(t)\coloneqq M(M(t,x),M(y,t)), \quad t \in I.
\end{equation*}
We first state the algebraic compatibility used in the construction.
It describes how the maps $\Phi_{x,y}$ interact with $M$, in particular showing that each $\Phi_{x,y}$ preserves the operation.

\begin{lem}\label{lem:Phi-identities}
For $x, y, z, w, r, s, t \in I$,
\begin{equation}\label{eq:Phi-four}
 M(\Phi_{x,y}(t),\Phi_{z,w}(s)) =\Phi_{M(x,z),M(y,w)}(M(t,s)),
\end{equation}
\begin{equation}\label{eq:Phi-left}
 M(r,\Phi_{x,y}(t)) =\Phi_{M(r,x),M(r,y)}(M(r,t)).
\end{equation}
In particular, if $z=x$ and $w=y$, then
\begin{equation}\label{eq:Phi-endomorphism}
 \Phi_{x,y}(M(a,b))=M(\Phi_{x,y}(a),\Phi_{x,y}(b)),
 \qquad a,b\in I.
\end{equation}
\end{lem}

\begin{proof}
Three applications of bisymmetry give
\begin{align*}
M(\Phi_{x,y}(t),\Phi_{z,w}(s)) &= M\left(M(M(t,x),M(y,t)),M(M(s,z),M(w,s))\right)\\
&=M\left(M(M(t,x),M(s,z)),M(M(y,t),M(w,s))\right)\\
&=M\left(M(M(t,s),M(x,z)),M(M(y,w),M(t,s))\right),
\end{align*}
which is \eqref{eq:Phi-four}.
Taking $x = y = t = r$ in \eqref{eq:Phi-four}, using $\Phi_{r,r}(r)=r$, and renaming $(z,w,s)$ as $(x,y,t)$ gives \eqref{eq:Phi-left}.
Taking $(z,w)=(x,y)$ in \eqref{eq:Phi-four} and using reflexivity gives \eqref{eq:Phi-endomorphism}.
\end{proof}

The next lemma is the interval case of Tarski's fixed-point theorem~\cite[Theorem~1]{Tarski1955}.
We include its short order-theoretic proof; increasing here means nondecreasing, and no continuity is assumed.

\begin{lem}\label{lem:monotone-fixed-point}
Let $h:[a,b]\to[a,b]$ be increasing.
Then $h$ has a fixed point.
Its least fixed point is $\inf\{t\in[a,b]:h(t)\leq t\}$, and its greatest fixed point is $ \sup\{t\in[a,b]:h(t)\geq t\}$.
\end{lem}

\begin{proof}
Let $E \coloneqq \{t:h(t)\leq t\}$ and $c \coloneqq \inf E$.
The set $E$ is nonempty since $b\in E$.
For every $v\in E$, monotonicity of $h$ gives $h(c)\leq h(v)\leq v$, and hence $h(c)\leq c$.
Therefore $h(h(c))\leq h(c)$, so $h(c)\in E$ and $c\leq h(c)$.
Thus $h(c)=c$.
Every fixed point belongs to $E$ and is at least $c$.
The greatest fixed point follows by reversing the order.
\end{proof}

\begin{thm}\label{thm:unique-fixed-point}
For every $x,y\in I$, the map $\Phi_{x,y}$ has exactly one fixed point, and it lies between $x$ and $y$.
\end{thm}

\begin{proof}
If $x=y$, then $x$ is a fixed point.
If $x<y$, the strict mean property gives that $x < \Phi_{x,y}(x) < \Phi_{x,y}(y) < y$.
Thus $h\coloneqq\Phi_{x,y}|_{[x,y]}$ maps $[x,y]$ into itself and is increasing, so Lemma~\ref{lem:monotone-fixed-point} gives a fixed point in $[x,y]$, and the strict endpoint inequalities place it in $(x,y)$.
The case $y<x$ is dealt with in the same manner.

For uniqueness, fix $x,y$ and write $h=\Phi_{x,y}$.
Suppose that $a$ and $b$ are fixed points of $h$, and let $c \coloneqq M(a,b)$.
By \eqref{eq:Phi-endomorphism},
\[ h(c)=M(h(a),h(b))=M(a,b)=c.\]
Applying \eqref{eq:Phi-left} first with $(r,t)=(a,b)$ and then with $(r,t)=(c,c)$ gives
\begin{equation}\label{eq:unique-a}
 c = M(a,b) = M(a, h(b)) = \Phi_{M(a,x),M(a,y)}(c),
\end{equation}
\begin{equation}\label{eq:unique-c}
 c = M(c,c) = M(c, h(c)) = \Phi_{M(c,x),M(c,y)}(c).
\end{equation}

For this fixed $c$, the function $r \mapsto\Phi_{M(r,x),M(r,y)}(c)$ is strictly increasing.
Therefore, equality of the right-hand sides of \eqref{eq:unique-a} and \eqref{eq:unique-c} gives $a=c$.
Hence $M(a,b)=a=M(a,a)$, and strict increase in the second variable yields $b=a$.
\end{proof}

Let  $m(x,y)$ be the unique fixed point of $\Phi_{x,y}$.

Our next proof builds substantially on Acz\'el's argument for non-symmetric means~\cite[pp.~396--398]{Aczel1948}.
In particular, the compatibility identity~\eqref{eq:intertwining-binary} below is equation~\textup{(15)} in his proof~\cite[p.~398]{Aczel1948}.

\begin{prop}\label{prop:mu-properties}
The operation $m:I^2\to I$ is reflexive, symmetric, strictly increasing in each variable and bisymmetric.
Moreover,
\begin{equation}\label{eq:intertwining-binary}
 M(m(x,y),m(z,w)) =m(M(x,z),M(y,w))
\end{equation}
for all $x,y,z,w\in I$.
\end{prop}

\begin{proof}
Reflexivity follows from $\Phi_{x,x}(x)=x$.
Bisymmetry of $M$ gives that
\[ \Phi_{x,y}(t) =M(M(t,y),M(x,t)) = \Phi_{y,x}(t),\]
so uniqueness of the fixed point implies symmetry of $m$.

In order to prove strict increase, suppose that $x\leq x'$ and $y\leq y'$, with at least one strict inequality.
Let $b \coloneqq m(x',y')$.
By Theorem \ref{thm:unique-fixed-point}, $b$ lies between $x^{\prime}$ and $y^{\prime}$.
Strict parameter monotonicity gives
\begin{equation}\label{eq:mu-strict}
 \Phi_{x,y}(b)<\Phi_{x',y'}(b)=b.
\end{equation}
Let $a \coloneqq \min\{x,y\}$.
Then $a<b$.
Indeed, $a \le \min\{x',y'\} \le b$, and $a=b$ would give $b=\min\{x',y'\}$, hence $x'=y'=b$ by the strict mean property and then $x=y=b$, contradicting the assumption.

Furthermore, $\Phi_{x,y}(a)\geq a$.
Hence $\Phi_{x,y}$ maps $[a,b]$ into itself.
Lemma~\ref{lem:monotone-fixed-point} supplies a fixed point in this interval, which must be $m(x,y)$.
Equation~\eqref{eq:mu-strict} shows that $b$ is not that fixed point, and therefore
$ m(x,y)<m(x',y')$.

\begin{sloppypar}
Equation~\eqref{eq:Phi-four} shows that the left-hand side of \eqref{eq:intertwining-binary} is fixed by $\Phi_{M(x,z),M(y,w)}$.
Uniqueness proves \eqref{eq:intertwining-binary}.
\end{sloppypar}

For bisymmetry of $m$, put
\[
 \begin{gathered}
 A=m(x,y),\quad B=m(z,w),\\
 C=m(x,z),\quad D=m(y,w),\quad T=m(A,B).
 \end{gathered}
\]
Repeated use of \eqref{eq:intertwining-binary} gives
\begin{align*}
 \Phi_{C,D}(T) &=M(M(T,C),M(D,T))\\
 &=M\left(m(M(A,x),M(B,z)), m(M(y,A),M(w,B))\right)\\
 &=m\left(M(M(A,x),M(y,A)),  M(M(B,z),M(w,B))\right)\\
 &=m(\Phi_{x,y}(A),\Phi_{z,w}(B))
 =m(A,B)=T.
\end{align*}
Thus $T$ is the unique fixed point of $\Phi_{C,D}$, so
\[ m(m(x,y),m(z,w)) =m(m(x,z),m(y,w)). \qedhere\]
\end{proof}

We now use the symmetric automatic-continuity theorem.

\begin{thm}[Burai, Kiss and Szokol~\cite{BuraiKissSzokol2021}]\label{thm:symmetric-input}
Let $N:I^2\to I$ be reflexive, symmetric, bisymmetric and strictly increasing in each variable.
Then there exist an interval $J\subset\R$ and a continuous strictly increasing bijection $\varphi:I\to J$ such that
\begin{equation*}
 \varphi(N(x,y))=\frac{\varphi(x)+\varphi(y)}2.
\end{equation*}
\end{thm}

\begin{proof}[Proof of Theorem~\ref{thm:main}]
Apply Theorem~\ref{thm:symmetric-input} to the operation $m$ from Proposition~\ref{prop:mu-properties}.
Let $J \coloneqq \varphi(I)$ and define
\[ G(s,t)\coloneqq \varphi\left(M(\varphi^{-1}(s),\varphi^{-1}(t))\right),  \qquad s,t\in J.\]
This operation is strictly increasing in each variable and satisfies $G(t,t)=t$.
We first explain how \eqref{eq:intertwining-binary} becomes Jensen's midpoint equation.
For $s,t,u,v\in J$, let
\[ x=\varphi^{-1}(s),\quad y=\varphi^{-1}(u),\quad z=\varphi^{-1}(t),\quad w=\varphi^{-1}(v).\]
Since $\varphi(m(x,y))=(s+u)/2$ and $\varphi(m(z,w))=(t+v)/2$, applying $\varphi$ to the two sides of \eqref{eq:intertwining-binary} gives, respectively,
\[ G\left(\frac{s+u}{2},\frac{t+v}{2}\right) \quad\hbox{and}\quad \frac{G(s,t)+G(u,v)}2.\]
Consequently,
\begin{equation*}
 G\left(\frac{s+u}{2},\frac{t+v}{2}\right)= \frac{G(s,t)+G(u,v)}{2},  \qquad s,t,u,v\in J.
\end{equation*}

We use the classical fact that an increasing solution of Jensen's midpoint equation on a real interval is affine; see~\cite[Section~13.2]{Kuczma2009}, and compare the Jensen-equation argument in~\cite{Maksa2002}.
For completeness, we give the elementary argument needed here.
Suppose $g:J\to\R$ is increasing and
\[ g\left(\frac{a+b}{2}\right)=\frac{g(a)+g(b)}{2},  \qquad a,b\in J.\]
Fix $a<b$ in $J$.
Induction on $k$ gives
\[ g\left((1-r)a+rb\right)=(1-r)g(a)+rg(b) \qquad\left(r=\frac{j}{2^k},\ 0 \leq j \leq 2^k\right).\]
For arbitrary $0<r<1$, choose dyadic numbers $r_k\leq r\leq s_k$ with $\lim_{k \to \infty} r_k = \lim_{k \to \infty} s_k = r$.
Monotonicity bounds $g((1-r)a+rb)$ between
\[ (1-r_k)g(a)+r_kg(b) \quad\hbox{and}\quad (1-s_k)g(a)+s_kg(b).\]
Both bounds tend to $(1-r)g(a)+rg(b)$.
Thus $g$ agrees with the affine interpolant on $[a,b]$.
The same identity on any larger compact subinterval shows that this affine formula holds throughout $J$.
In particular, continuity is a consequence of midpoint preservation and monotonicity, not an assumption in this step.

Setting $v=t$ in the midpoint equation for $G$ makes each first-variable section an increasing midpoint-preserving function.
The preceding argument makes it affine; setting $u=s$ gives the same conclusion for every second-variable section.
In order to see the resulting form of $G$ explicitly, fix $a<b$ in $J$.
Affinity in the first variable gives
\[ G(s,t)=\frac{b-s}{b-a}G(a,t)+\frac{s-a}{b-a}G(b,t),  \qquad s,t\in J.\]
Both $G(a,t)$ and $G(b,t)$ are affine functions of $t$.
Therefore expanding yields
\begin{equation*}
 G(s,t)=Ast+Bs+Ct+D
\end{equation*}
with constant real coefficients.
Reflexivity now says
\[ At^2+(B+C)t+D=t, \qquad t\in J.\]
Since $J$ is non-degenerate, this polynomial identity implies $A=D=0$ and $B+C=1$.
Strict increase gives $B,C>0$, so with $p=B$ we obtain
\[ G(s,t)=ps+(1-p)t,\qquad 0<p<1.\]
Conjugating back proves \eqref{eq:main-representation}.
The continuous strictly increasing bijection $\varphi:I\to J$ has a continuous inverse, so this representation also proves continuity of $M$.

For uniqueness, suppose that $M$ has another representation with weight $q\in(0,1)$ and continuous increasing generator $\psi$.
The auxiliary mean $m$ is defined by $M$ itself, independently of a representation.
In any representation with weight $r\in(0,1)$ and generator $\theta$, put
\[ z=\theta^{-1}\!\left(\frac{\theta(x)+\theta(y)}2\right).
\]
Substitution into the definition of $\Phi_{x,y}$ gives
\[  \theta(\Phi_{x,y}(z))=(r^2+(1-r)^2)\theta(z) +r(1-r)(\theta(x)+\theta(y)) =\theta(z).\]
Thus $z$ is the unique fixed point, and
\[ \theta(m(x,y))=\frac{\theta(x)+\theta(y)}2.\]
In particular, both $\varphi$ and $\psi$ transform the same mean $m$ into the arithmetic midpoint.
Their transition map $h=\psi\circ\varphi^{-1}$ therefore satisfies
\[ h\left(\frac{s+t}{2}\right)=\frac{h(s)+h(t)}{2},  \qquad s,t\in J.\]
It is increasing, so the midpoint argument just proved gives $h(s)=as+b$ on $J$, with $a>0$.
Hence $\psi=a\varphi+b$ throughout $I$.
The weighted mean is unchanged by this affine transformation, so
\[
 (p-q)(s-t)=0, \qquad s,t\in J.
\]
\nopagebreak[4]
As $J$ is non-degenerate, $p=q$.
\end{proof}

The construction also identifies the auxiliary fixed point explicitly.

\begin{cor}\label{cor:explicit-fixed-point}
With the notation of Theorem~\ref{thm:main},
\begin{equation}\label{eq:mu-explicit}
 m(x,y) = \varphi^{-1}\!\left(\frac{\varphi(x)+\varphi(y)}2\right).
\end{equation}
Moreover,
\begin{equation}\label{eq:Phi-affine}
 \varphi(\Phi_{x,y}(t)) = (p^2+(1-p)^2)\varphi(t)+p(1-p)(\varphi(x)+\varphi(y)),
\end{equation}
and hence
\begin{equation}\label{eq:Phi-sign}
 \varphi(\Phi_{x,y}(t))-\varphi(t) = p(1-p)\left(\varphi(x)+\varphi(y)-2\varphi(t)\right).
\end{equation}
\end{cor}

Thus, the unique fixed point of $\Phi_{x,y}$ is the quasi-arithmetic midpoint \eqref{eq:mu-explicit}, and $\Phi_{x,y}(t) > t$ for $t < m(x,y)$, $\Phi_{x,y}(t) < t$ for $t > m(x,y)$.
The map $\Phi_{x,y}$ pushes every point towards $m(x,y)$.


\begin{proof}
The representation in Theorem~\ref{thm:main} gives
\[ \varphi(M(t,x))=p\varphi(t)+(1-p)\varphi(x), \qquad \varphi(M(y,t))=p\varphi(y)+(1-p)\varphi(t).\]
Applying the same representation to the outer occurrence of $M$, we obtain
\begin{align*}
 \varphi(\Phi_{x,y}(t))
 &=p\varphi(M(t,x))+(1-p)\varphi(M(y,t))\\
 &=p\left(p\varphi(t)+(1-p)\varphi(x)\right) +(1-p)\left(p\varphi(y)+(1-p)\varphi(t)\right)\\
 &=(p^2+(1-p)^2)\varphi(t)
       +p(1-p)(\varphi(x)+\varphi(y)),
\end{align*}
which proves \eqref{eq:Phi-affine}.
Since $p^2+(1-p)^2-1=-2p(1-p)$, subtracting $\varphi(t)$ gives \eqref{eq:Phi-sign}.
The map $\varphi$ is injective and $p(1-p)>0$, so
\[\Phi_{x,y}(t)=t \quad\Longleftrightarrow\quad \varphi(t)=\frac{\varphi(x)+\varphi(y)}2.\]
The midpoint on the right belongs to the interval $J=\varphi(I)$.
Therefore, it has a unique inverse image in $I$, proving \eqref{eq:mu-explicit}.
Finally, strict increase of $\varphi$ shows that the sign of
$\Phi_{x,y}(t)-t$ agrees with the sign of
$\varphi(x)+\varphi(y)-2\varphi(t)$.
This sign is positive for $t<m(x,y)$, zero at $m(x,y)$, and negative for $t>m(x,y)$.
\end{proof}

\begin{rem}
Acz\'el~\cite[pp.~396--397]{Aczel1948} constructs the symmetric mean by the compound-mean iteration
\[ a_{n+1} = \min\{M(a_n,b_n),M(b_n,a_n)\}, b_{n+1} = \max\{M(a_n,b_n),M(b_n,a_n)\}, \]
starting from $a_0=\min\{x,y\}$ and $b_0=\max\{x,y\}$.
Continuity of $M$ ensures that the two limits agree, but this argument is unavailable before continuity has been established.

Without assuming continuity, one may restart the iteration from the limiting pair and continue through countable ordinals, taking monotone limits at limit stages.
This argument appeared in \cite[Lemma 3]{Pasteczka2019}.
The strict mean property and the countability of $\Q$ ensure that the endpoints meet at some countable stage.
This transfinite construction uses only reflexivity and strict increase to define a symmetric mean.
As in Acz\'el's approach, its purpose is to introduce such a mean as a first step towards proving Theorem~\ref{thm:main}, here without assuming continuity.
The representation \eqref{eq:main-representation}, once it is established, identifies this mean with $m(x,y)$, which is the unique fixed point of $\Phi_{x,y}$.
It also shows that the ordinary iteration already converges to this common value.
\end{rem}




\section{The higher-arity consequence}\label{sec:nary}

In the symmetric case, Kiss and Shulman~\cite{KissShulman} obtain automatic continuity for operations of arbitrary finite arity by constructing an averaging scale on a dense set.
Theorem~\ref{thm:nary-main} extends this higher-arity characterization to operations without a symmetry assumption.

We follow the binary-retract reduction of Maksa~\cite[Section~3]{Maksa2002}.
Theorem~\ref{thm:main} supplies its binary input without a continuity assumption and allows the retract to be used even when it is not symmetric.
The density argument below likewise obtains the needed Jensen conclusion from monotonicity alone.

Let $n\geq2$.
An operation $F:I^n\to I$ is reflexive if $F(x,\ldots,x)=x$ for every $x \in I$, and bisymmetric if, for every $n\times n$ array $(x_{ij})$ in $I$,
\begin{align}\label{eq:nary-bisymmetry}
&F\left(F(x_{11},\ldots,x_{1n}),\ldots,
        F(x_{n1},\ldots,x_{nn})\right)\notag\\
&\qquad=
F\left(F(x_{11}, x_{21}, \ldots,x_{n1}),\ldots,
        F(x_{1n}, x_{2n}, \ldots,x_{nn})\right),
\end{align}
and strictly increasing in each variable if every one‑variable section of $F$ is strictly increasing.

Define the {\it binary retract} $ B(x,y)\coloneqq F(x,y,\ldots,y)$.

\begin{lem}
\label{lem:binary-retract}
If $F$ is reflexive, bisymmetric and strictly increasing in each variable, then $B$ is reflexive, bisymmetric and strictly increasing in each variable.
Moreover,
\begin{equation}\label{eq:nary-intertwining}
F(B(x_1,y_1),\ldots,B(x_n,y_n)) = B(F(x_1,\ldots,x_n),F(y_1,\ldots,y_n)).
\end{equation}
\end{lem}

\begin{proof}
Reflexivity and strict increase are immediate.
Apply \eqref{eq:nary-bisymmetry} to the array
\[
\begin{pmatrix}
x&y&\cdots&y\\
u&v&\cdots&v\\
\vdots&\vdots&&\vdots\\
u&v&\cdots&v
\end{pmatrix}.
\]
Evaluation by rows and columns gives
\[ B(B(x,y),B(u,v))=B(B(x,u),B(y,v)),\]
so $B$ is bisymmetric.
For \eqref{eq:nary-intertwining}, use the array whose $i$th row is $(x_i,y_i,\ldots,y_i)$.
Applying \eqref{eq:nary-bisymmetry} to this array gives \eqref{eq:nary-intertwining}.
\end{proof}

The following lemma uses a subdivision and density argument similar to the one employed by Kiss and Shulman~\cite{KissShulman}.

\begin{lem}\label{lem:fixed-weight-affine}
Let $J$ be a non-degenerate real interval, let $0<\rho<1$, let $m \ge 1$, and let $H:J^m\to\R$ be strictly increasing in each variable.
Suppose
\begin{equation}\label{eq:vector-Jensen}
 H(\rho\bm{x}+(1-\rho)\bm{y}) = \rho H(\bm{x})+(1-\rho)H(\bm{y})
\end{equation}
for all $\bm{x},\bm{y}\in J^m$.
Then
\[ H(u_1,\ldots,u_m) = c+\lambda_1 u_1+\cdots+\lambda_m u_m, \qquad u_1, \dots, u_m \in J, \]
for some $c\in\R$ and positive numbers $\lambda_1,\ldots,\lambda_m$.
\end{lem}

\begin{proof}
Fix all variables except the $i$th one, and denote the resulting one-variable section by $g:J\to\R$.
Then $g$ is strictly increasing, and \eqref{eq:vector-Jensen} gives
\[ g(\rho u+(1-\rho)v)=\rho g(u)+(1-\rho)g(v),  \qquad u,v\in J.\]
We show directly that $g$ is affine, without assuming continuity.
Fix $a<b$ in $J$ and let $\ell$ be the affine function agreeing with $g$ at $a,b$, specifically,
\[\ell(t) = \frac{b-t}{b-a}\,g(a) + \frac{t-a}{b-a}\,g(b),  \qquad t \in J.\]
Subdivide $[a,b]$ by inserting $\rho a+(1-\rho)b$, and repeat this subdivision in every resulting interval.
The fixed-weight identity shows inductively that $g=\ell$ at every subdivision point.
After $k$ rounds, the largest interval has length at most $(b-a)\max\{\rho,1-\rho\}^{k}$, which tends to zero as $k \to \infty$.
For any $t\in[a,b]$, take subdivision endpoints $a_k\leq t\leq b_k$ with $a_k,b_k\to t, k \to \infty$.
Monotonicity gives that
\[ \ell(a_k)=g(a_k)\leq g(t)\leq g(b_k)=\ell(b_k).\]
By $k \to \infty$, we have $g(t)=\ell(t)$.
The affine formulas on overlapping compact subintervals agree, so $g$ is affine on $J$.
Thus every one-variable section of $H$ is affine.

Separate affinity implies that $H$ is multi-affine:
\begin{equation}\label{eq:multiaffine}
 H(\bm{u})=\sum_{A\subseteq\{1,\ldots,m\}}c_A\prod_{i\in A}u_i.
\end{equation}
Here the $c_A$ are real constants, and the term for $A=\emptyset$ is $c_\emptyset$, with the empty product equal to $1$.
The case $m=1$ is immediate.
For $m\geq2$, fix $a<b$ in $J$ and write $\bm{u}'=(u_1,\ldots,u_{m-1})$.
Affinity in the last variable gives
\[ H(\bm{u}',t)=\frac{b-t}{b-a}H(\bm{u}',a)             +\frac{t-a}{b-a}H(\bm{u}',b),  \qquad t\in J.\]
The two sections on the right are separately affine in $m-1$ variables, so induction yields \eqref{eq:multiaffine}.

Finally, substitute \eqref{eq:multiaffine} into \eqref{eq:vector-Jensen}.
Since $J$ has nonempty interior, this is an identity between polynomials on a nonempty open set, so their coefficients agree.
If $|A|\geq2$, split $A=A_1 \cup A_2, A_1 \cap A_2 = \emptyset, A_1 \ne \emptyset, A_2 \ne \emptyset$.
In the expansion of $H(\rho\bm{x}+(1-\rho)\bm{y})$, the coefficient of the mixed monomial $\displaystyle \prod_{i\in A_1}x_i\prod_{j\in A_2}y_j$ is $c_A\rho^{|A_1|}(1-\rho)^{|A_2|}$,
whereas no such monomial occurs in the expansion of $\rho H(\bm{x})+(1-\rho)H(\bm{y})$.
Hence $c_A=0$ if $|A| \geq 2$.
Only the constant and linear terms remain, and strict increase makes every linear coefficient positive.
\end{proof}

\begin{proof}[Proof of Theorem~\ref{thm:nary-main}]
First suppose that $F$ is reflexive, bisymmetric and strictly increasing in each variable.
Apply Theorem~\ref{thm:main} to the binary retract $B$.
There are $\rho\in(0,1)$, an interval $J$ and a continuous increasing bijection $\varphi:I\to J$ such that
\[ \varphi(B(x,y))=\rho\varphi(x)+(1-\rho)\varphi(y).\]
For $\bm{u}=(u_1,\ldots,u_n)\in J^n$, define
\begin{equation}\label{eq:H-definition}
 H(\bm{u})\coloneqq \varphi\left(F(\varphi^{-1}(u_1),\ldots,\varphi^{-1}(u_n))\right).
\end{equation}
In the coordinates~\eqref{eq:H-definition}, equation~\eqref{eq:nary-intertwining} becomes
\[ H(\rho\bm{u}+(1-\rho)\bm{v}) =\rho H(\bm{u})+(1-\rho)H(\bm{v}), \quad \bm{u}, \bm{v} \in J^n.\]
Therefore Lemma~\ref{lem:fixed-weight-affine} yields that
\[ H(u_1,\ldots,u_n)=c+\lambda_1u_1+\cdots+\lambda_nu_n,  \qquad u_1, \dots, u_n \in J, \]
for some $c \in \mathbb{R}$ and $\lambda_1, \dots, \lambda_n > 0$.
Reflexivity says $H(t,\ldots,t)=t$ on the non-degenerate interval $J$, hence, $c=0$ and $\lambda_1+\cdots+\lambda_n=1$.
Conjugating back proves \eqref{eq:nary-form} and continuity.

For uniqueness, a second representation restricts on the binary retract to a second representation of $B$.
The uniqueness part of Theorem~\ref{thm:main} identifies its first weight and shows that the two generators differ by a positive affine change.
After making this change, variation of one coordinate at a time in the two linear forms gives equality of all the weights.

Conversely, suppose \eqref{eq:nary-form} holds with positive weights summing to $1$ and a continuous strictly increasing generator.
The weighted sum lies in the interval $\varphi(I)$, so the operation is well-defined, reflexive and strictly increasing in each variable.
For an arbitrary array $(x_{ij})$, applying $\varphi$ to either iterated evaluation gives
\[ \sum_{i=1}^n\sum_{j=1}^n\lambda_i\lambda_j\varphi(x_{ij}),
\]
where the two sums are interchanged for column evaluation.
The expressions agree, and injectivity of $\varphi$ proves bisymmetry.
\end{proof}

\section{Strongly decomposable families without regularity}\label{sec:strong-decomposability}

As an application of our results, we remove the continuity assumption from Marichal's characterization \cite[Theorem~2.3]{Marichal2000}.

A family $(M_n)_{n\geq1}$ of operations $M_n:I^n\to I$ is \emph{strongly decomposable} if replacing the coordinates indexed by any nonempty subset $K=\{i_1<\cdots<i_k\}\subseteq\{1,\ldots,n\}$ by repeated copies of $M_k(x_{i_1},\ldots,x_{i_k})$ leaves $M_n(x_1,\ldots,x_n)$ unchanged.

We first extend \cite[Proposition~2.1]{Marichal2000} to arbitrary arity.

\begin{lem}\label{lem:strong-decomposability-bisymmetry}
Let $(M_k)_{k\geq1}$ be a strongly decomposable family of reflexive operations $M_k:I^k\to I$.
Then $M_n$ is bisymmetric for every $n\geq2$.
\end{lem}

\begin{proof}
Fix $n\geq2$ and an $n\times n$ array $(x_{ij})$ in $I$.
Let $\bm{x}$ be the $n^2$-tuple obtained by reading the array row by row, and put
\[ r_i=M_n(x_{i1},\ldots,x_{in}),\qquad c_j=M_n(x_{1j},\ldots,x_{nj}).\]
Then
\begin{align*}
 M_n(r_1,\ldots,r_n)
 &=M_{n^2}(\bm{x})\\
 &=M_{n^2}(\underbrace{c_1,\ldots,c_n;\ \ldots;\ c_1,\ldots,c_n}_{n\text{ identical rows}})\\
 &=M_n(c_1,\ldots,c_n).
\end{align*}
The first and last equalities follow from \cite[Lemma~2.1, equation~(6)]{Marichal2000} and reflexivity.
The middle equality follows by applying strong decomposability successively to the disjoint columns, replacing each entry in column $j$ by $c_j$.
Thus the row and column evaluations agree.
\end{proof}

\begin{thm}\label{thm:strong-decomposability}
Let $I$ be a non-degenerate real interval.
A family $(M_n)_{n\geq1}$ of operations $M_n:I^n\to I$ is strongly decomposable, reflexive and strictly increasing in each variable if and only if there exists a continuous strictly increasing function $\varphi:I\to\R$ such that
\[M_n(x_1,\ldots,x_n) =\varphi^{-1}\!\left(\frac1n\sum_{i=1}^n\varphi(x_i)\right), \qquad n\geq1,\quad x_1,\ldots,x_n\in I.\]
In particular, every $M_n$ is continuous and symmetric.
\end{thm}

The proof strategy is the same as in the proof of \cite[Theorem~2.3]{Marichal2000}.

\begin{proof}
The converse is immediate.
Suppose that the family has the stated properties.
By Lemma~\ref{lem:strong-decomposability-bisymmetry}, $M_2$ is bisymmetric.
Theorem~\ref{thm:main}, which removes continuity from \cite[Theorem~2.1]{Marichal2000}, therefore gives
\[ M_2(x,y)=\varphi^{-1}\!\left(p\varphi(x)+(1-p)\varphi(y)\right),
 \qquad 0<p<1,\]
with $\varphi$ continuous and strictly increasing.

Lemma~\ref{lem:strong-decomposability-bisymmetry} and Theorem~\ref{thm:nary-main} also show that $M_n$ is continuous for every $n\geq2$.
In particular, continuity of $M_3$ follows from the ternary case of Theorem~\ref{thm:nary-main}.

Let $J \coloneqq \varphi(I)$ and
\[ F_k(z_1,\ldots,z_k) \coloneqq \varphi\!\left(M_k(\varphi^{-1}(z_1),\ldots,\varphi^{-1}(z_k))\right), \qquad k=2,3.\]
Then $F_2(z_1,z_2)=pz_1+(1-p)z_2$, and $F_3$ is continuous.
Consequently, the matrix iteration in the proof of \cite[Theorem~2.3]{Marichal2000}, using \cite[Lemma~2.2]{Marichal2000}, applies without change and gives
\[ F_3(z_1,z_2,z_3) =\frac{p^2z_1+p(1-p)z_2+(1-p)^2z_3}{p^2+p(1-p)+(1-p)^2}.\]
Strong decomposability on the first and third coordinates now gives, as in that proof,
\[ p(1-p)(1-2p)(z_3-z_1)=0, \qquad z_1,z_3\in J.\]
Hence $p=1/2$, and $M_2$ is symmetric, and by iteration
(see also \cite[Lemma~2.3]{Marichal2000}), every $M_n$ is symmetric.
Since their continuity has already been established (and $M_1$ is the identity by reflexivity), the Kolmogorov--Nagumo characterization \cite[Theorem~2.2]{Marichal2000} completes the proof.
\end{proof}

\section{Further questions}\label{sec:questions}

\subsection{Weighted compatible families}

The Kolmogorov--Nagumo--de~Finetti characterization concerns a compatible family of operations, whereas Theorem~\ref{thm:nary-main} treats one fixed arity at a time.
A weighted counterpart should identify a single continuous strictly increasing generator $\varphi$ for a family $(M_n)_{n\geq2}$, with representations
\begin{equation}\label{eq:weighted-family}
 M_n(x_1,\ldots,x_n) =\varphi^{-1}\!\left(\sum_{i=1}^n\lambda_{n,i}\varphi(x_i)\right),
 \qquad x_1, \dots, x_n \in I,
\end{equation}
for $\lambda_{n,i}>0, 1 \le i \le n, \quad \sum_{i=1}^n\lambda_{n,i}=1$.
The cross-arity assumptions matter: applying Theorem~\ref{thm:nary-main} separately does not by itself identify the generators for different $n$.
Classical axiomatizations without symmetry, including~\cite{Marichal2000}, provide a natural point of comparison.

One sufficient algebraic compatibility condition is rectangular bisymmetry: for every $m,n\geq2$ and every $m\times n$ array $(x_{ij})$,
\begin{align}\label{eq:rectangular-family}
 &M_m\left(M_n(x_{11},x_{12},\ldots,x_{1n}),\ldots,
                 M_n(x_{m1},x_{m2},\ldots,x_{mn})\right)\notag\\
 &\qquad =M_n\left(M_m(x_{11},x_{21},\ldots,x_{m1}),\ldots,
                 M_m(x_{1n},x_{2n},\ldots,x_{mn})\right).
\end{align}
Indeed, assume that the $M_n$ are reflexive and strictly increasing in each variable and satisfy~\eqref{eq:rectangular-family}.
The case $m=n=2$ and Theorem~\ref{thm:main} give a generator for $M_2$.
In those coordinates, the case $m=2$ of~\eqref{eq:rectangular-family} is the fixed-weight Jensen identity for each transformed $M_n$.
Lemma~\ref{lem:fixed-weight-affine} and reflexivity therefore yield~\eqref{eq:weighted-family} with that same generator, without any continuity assumption.
Conversely, families of the form~\eqref{eq:weighted-family} satisfy~\eqref{eq:rectangular-family}.
Thus this stronger compatibility condition already gives a complete characterization without regularity.

\begin{question}\label{ques:weighted-family}
What weighted block-replacement axioms can take the place of the rectangular identities~\eqref{eq:rectangular-family} in a Kolmogorov--Nagumo--de~Finetti type characterization without continuity assumptions?
Which systems of positive weights are compatible with such axioms?
\end{question}

Here a weighted block replacement should use the relative weights of that block and assign their sum to the resulting block mean.
If arbitrary probability vectors are used as indices, the target representation is
\[
 M_{\bm{w}}(\bm{x})
 =\varphi^{-1}\!\left(\sum_{i=1}^{n} w_i\varphi(x_i)\right),
 \qquad x_1, \dots, x_n \in I,
\]
where $\bm{w} = (w_1, \dots, w_n)$ such that $ w_i>0, 1 \le i \le n, \ \sum_i w_i=1$.
The question is to identify suitable replacement axioms and the extent to which they can be weakened while still forcing a common generator.

\subsection{Nonstrict monotonicity and strict means}

The strict increase assumption in Theorem~\ref{thm:main} cannot be replaced by nondecrease, even if symmetry and the strict mean property are assumed.

\begin{exa}\label{exa:strict-mean}
On $I=[-1,1]$, define
\[ N(x,y) \coloneqq
 \begin{cases}
  (x+y)/2,&x,y<0\ \text{or}\ x,y\geq0,\\[2pt]
  \min\{x,y\}/2,&\min\{x,y\}<0\leq\max\{x,y\}.
 \end{cases}
\]
Then $N$ is a reflexive, symmetric, bisymmetric strict mean, nondecreasing in each variable, but it is neither strictly increasing in each variable nor continuous.
\end{exa}

\begin{proof}
Symmetry and reflexivity are immediate.
Assume $x < y$.
If $x, y \ge 0$ or $x, y < 0$, the strict mean inequality $x < N(x,y) < y$ follows from the inequality $x < \frac{x+y}{2} < y$.
If $x < 0 \le y$, it follows from the inequality $x<x/2<y$.

For fixed $y<0$, the section $x\mapsto N(x,y)$ equals $(\min\{x,0\}+y)/2$, which is not strictly increasing on $[0,1]$.
For fixed $y\geq0$, it equals $x/2$ on $x<0$ and $(x+y)/2$ on $x\geq0$, which is discontinuous at $x = 0$ if $y > 0$.
These sections are nondecreasing, and symmetry gives the same property in the other variable.

We show bisymmetry.
Let $t_{-} \coloneqq \min\{t,0\}$.
If $\min\{x,y\}<0$, then both cases of the definition give $N(x,y) = (x_-+y_-)/2 < 0$, so $N(x,y)_- = N(x,y)$.
If $x,y \ge 0$, then $N(x,y) \ge 0$ and $x_- = y_- = 0$.
Hence,
\begin{equation}\label{eq:N-negative-part}
N(x,y)_{-} = \frac{x_- + y_-}{2}, \qquad x,y \in I.
\end{equation}

Since $t<0$ if and only if $t_-<0$, \eqref{eq:N-negative-part} also shows that $N(x,y)<0$ if and only if $\min\{x,y\}<0$.

If $x,y,u,v$ are all nonnegative, then
\[
N(N(x,y),N(u,v)) = \frac{N(x,y) + N(u,v)}{2} = \frac{x + y + u + v}{4} = N(N(x,u),N(y,v)).
\]

If at least one of them is negative, then, by \eqref{eq:N-negative-part}, $N(N(x,y),N(u,v)) < 0$ and
\[ N(N(x,y),N(u,v)) = \frac{N(x,y)_{-} + N(u,v)_{-}}{2} = \frac{x_{-} + y_{-} + u_{-} + v_{-}}{4},\]  which is equal to $N(N(x,u),N(y,v))$.
Now bisymmetry follows.
\end{proof}

The continuous theory of nonstrict symmetric means is studied in~\cite{FodorMarichal1997}.
With continuity, strict mean property recovers strict increase.
Symmetry is not needed for this implication.
We now adopt their method in our concrete situation.

\begin{prop}\label{prop:continuous-strict-mean}
Let $I$ be a non-degenerate real interval and let $M:I^2\to I$ be continuous and bisymmetric.
Then $M$ is a strict mean if and only if it is reflexive and strictly increasing in each variable. \\
In particular, every continuous bisymmetric strict mean is a weighted quasi-arithmetic mean.
\end{prop}

\begin{proof}
Reflexivity and strict increase imply the strict mean property by the inequalities in the Introduction.

For the converse, assume that $M$ is a strict mean.
Then $M$ is reflexive, so it remains to prove strict increase in each variable. 
Fix $a\in I$ and set $g(t) \coloneqq M(a,t)$.
Reflexivity and bisymmetry give
\[ g(M(s,t)) = M(M(a,a), M(s,t)) = M(g(s),g(t)), \qquad s,t\in I.\]

Suppose that $g(x) = g(y) = c$ for some $x<y$.
Strict mean property gives that $a<g(t)<t$ for $t>a$, $t<g(t)<a$ for $t<a$, and $g(a)=a$.
Hence $g(x)=g(y)$ is impossible when $x \le a \le y$.
Hence, either $a<x<y$ or $x<y<a$.
Replacing $I$ by $-I$ and $M$ by $(s,t) \mapsto -M(-s,-t)$, which preserves all hypotheses, we may assume $a<x<y$.

By continuity, the number
\[z \coloneqq \min\{t\in[a,x]:g(t)=c\} \]
exists, 
and $z \ne a$ because $g(a) = a \ne c$.
Strict mean property applied to $g(x)$ and $g(z)$ gives $a < c < x$ and $c = g(z) < z$, so $a < c < z \le x < y$.
Since $M(a,y)=c<z<M(z,y)$, the intermediate value theorem supplies $t\in(a,z)$ with $M(t,y)=z$.
Consequently,
\[ c=g(z)=g(M(t,y))=M(g(t),g(y))=M(g(t),c).\]
The strict mean property forces $g(t)=c$, contradicting the definition of $z$.
Thus $g$ is injective.
By this and the continuity of $g$,
$g$ is strictly monotone.
Since $g(a)=a$ and $g(t)$ lies strictly between $a$ and $t$ for $t\ne a$, $g$ must be strictly increasing.

Applying the same argument to the opposite operation $(x,y)\mapsto M(y,x)$ gives strict increase in the first variable as well.
The representation now follows from Theorem~\ref{thm:main}.
\end{proof}

Therefore, Example~\ref{exa:strict-mean} shows that this equivalence fails without continuity, even in the symmetric, nondecreasing class.
It motivates the corresponding classification problem in the regularity-free setting.

\begin{question}\label{ques:nonstrict}
How can we characterize the symmetric, bisymmetric strict means $M:I^2\to I$ on a non-degenerate real interval which are nondecreasing in each variable, without assuming continuity?
\end{question}

More generally, this study can be extended to describe reflexive, bisymmetric operations on a real interval which are nondecreasing in each variable, without assuming continuity or symmetry or the strict mean property.

\vspace{1pc}

\noindent{\bf Use of generative AI.} \ The uniqueness part of the proof of Theorem  \ref{thm:unique-fixed-point} was suggested by ChatGPT 6 Astra.
The authors take full responsibility for the content and integrity of the publication.

\vspace{1pc}

\noindent{\bf Acknowledgements.}
The first author was supported by the Hungarian National Research, Development
and Innovation Office grants K146922, STARTING 150576. The second author has received funding from JSPS KAKENHI JP22K13928.

\end{document}